\documentclass[12pt,leqno]{article}
\usepackage{amssymb}
\usepackage[mathscr]{eucal}
\usepackage{amsmath,amssymb,latexsym,theorem,bbm}
\usepackage{amsmath,amssymb,latexsym,theorem}
\usepackage{color,url}
\usepackage{appendix}

\newcommand{\CC}{\mathbb{C}}

\newcommand{\NN}{\mathbb{N}}

\newcommand{\RR}{\mathbb{R}}

\newcommand{\ZZ}{\mathbb{Z}}

\newcommand{\cA}{{\mathcal A}}
\newcommand{\cB}{{\mathcal B}}

\newcommand{\cD}{{\mathcal D}}

\newcommand{\dd}{\mathrm{d}}
\newcommand{\ee}{\mathrm{e}}
\newcommand{\ii}{\mathrm{i}}

\newcommand{\EE}{\operatorname{\mathbb{E}}}

\newcommand{\PP}{\operatorname{\mathbb{P}}}

\newcommand{\vare}{\varepsilon}

\renewcommand{\leq}{\leqslant}
\renewcommand{\geq}{\geqslant}

\newcommand{\distre}{\stackrel{\cD}{=}}

\newcommand{\proofend}{\hfill\mbox{$\Box$}}

\numberwithin{equation}{section}

\theoremstyle{change} \theorembodyfont{\em}
\newtheorem{Lem}{Lemma.}[section]
\newtheorem{Thm}[Lem]{Theorem.}
\newtheorem{Pro}[Lem]{Proposition.}

\newtheorem{Def}[Lem]{Definition.}

\theorembodyfont{\rm}
\newtheorem{Rem}[Lem]{Remark.}
\newtheorem{Ex}[Lem]{Example.}

\begin{document}

\begin{center}
 {\bfseries\Large On calculation of the exact Lyapunov exponent }\\[1mm]
 {\bfseries\Large of affine Boole transformations}

\vspace*{3mm}

 {\sc\large
 M\'aty\'as $\text{Barczy}^{*}$}

\end{center}

\vskip0.2cm

\noindent
 * ELKH-SZTE Analysis and Applications Research Group,
   Bolyai Institute, University of Szeged,
   Aradi v\'ertan\'uk tere 1, H--6720 Szeged, Hungary.

\noindent e-mail: barczy@math.u-szeged.hu

%\noindent $\diamond$ Corresponding author.

{\renewcommand{\thefootnote}{}
\footnote{M\'aty\'as Barczy was supported by the project TKP2021-NVA-09.
Project no.\ TKP2021-NVA-09 has been implemented with the support
 provided by the Ministry of Innovation and Technology of Hungary from the National Research, Development and Innovation Fund,
 financed under the TKP2021-NVA funding scheme.}
}

\vspace*{-1cm}

\begin{center}
\textit{To the memory of Gabriella Vas}
\end{center}

\begin{abstract}
We show an elementary way for calculation of the exact Lyapunov exponent of affine Boole transformations using
 the interchangeability theorem for differentiation and integration due to Leibniz.
\end{abstract}

\section{Introduction and results}
\label{section_intro}

The theory and applications of Lyapunov exponents are important and extensively growing areas in mathematics and physics as well.
Roughly speaking, Lyapunov exponent of a dynamical system can be considered as a measure of the degree of sensitive dependence on the initial conditions of the dynamical system
 in question.
For several precise mathematical results on the relationships of Lyapunov exponents and sensitive dependence on initial conditions
 in the context of mappings of the unit interval into itself, see, e.g., Ko\c{c}ak and Palmer \cite{KocPal}.
Lyapunov exponents play an important role in calculating metric (also called measure-theoretic) entropy of mappings from
 a probability space to itself, in the theory of random matrices and random maps;
 in linear stochastic systems and stability theory; in random Schr\"odinger operators and wave propagation in random media;
 in nonlinear stochastic systems and stochastic flows on manifolds; and in chaos and phase transitions,
 see the lecture note of Arnold and Wihstutz \cite{ArnWih} and the very recent paper of Blackmore et al.\ \cite{BlaBalKycPry}.
The recent book of Barreira \cite{Bar} contains applications of Lyapunov exponents in hyperbolicity, ergodic theory, and multifractal analysis as well.

We point out the fact that, in general, one can find it difficult to obtain an explicit formula for Lyapunov exponent of a given dynamical system,
 and hence in practice numerical procedures are applied.
For example, very recently Miranda-Filho et al.\ \cite{Mir_et} have numerically approximated the so-called largest Lyapunov exponent
 for the Vicsek model that describes the dynamics of a flock of self-propelled particles.

In this paper we focus on affine Boole transformations, and hence throughout the paper
 we will recall the existing and corresponding literature only on the calculation of Lyapunov exponents of these mappings.
In the present paper, using Birkhoff's ergodic theorem for the Newton iteration for the mapping $x\mapsto (x-a)^2 + b^2$
 (where \ $a$ \ is a real number and \ $b$ \ is a positive real number)
 and the interchangeability theorem for differentiation and integration due to Leibniz in 1697
 (see, e.g., Engelsman \cite[page 9 and Chapter 3]{Eng}),
 we calculate the Lyapunov exponent of the corresponding affine Boole transformation
 \ $x\mapsto \frac{b}{2}\big( \frac{x-a}{b} - \frac{1}{\frac{x-a}{b}} \big) + a$, \ which is known to be \ $\ln(2)$.
\ The result itself is widely known, but our proof technique seems to be new, that's why we decided to write the present note.
The interchangeability theorem for differentiation and integration is also known as the method of differentiation under the integral sign
 or ''Feynman's trick'' (since Richard Feynman occasionally used this method for calculating complicated integrals, see Feynman \cite{Fey}).
For a good reference on the applications of the method of differentiation under the integral sign for calculating integrals, see Nahin \cite[Chapter 3]{Nah}.

Let \ $\ZZ_+$, \ $\NN$, \ $\RR$, \ $\RR_+$ \ and \ $\CC$ \ denote the
 set of non-negative integers, positive integers, real numbers, non-negative real numbers, and complex numbers, respectively.
The Borel \ $\sigma$-algebra on \ $\RR$ \ and \ $(0,\infty)$ \ is denoted by \ $\cB(\RR)$ \ and \ $\cB((0,\infty))$, \ respectively.
All the random variables are defined on a probability space \ $(\Omega,\cA,\PP)$.
\ Equality in distribution is denoted by \ $\distre$.
\ A random variable \ $\xi$ \ is said to have Cauchy distribution with parameter \ $(a,b)$, \ where \ $a\in\RR$ \ and \ $b>0$, \
 if \ $\xi$ \ has a density function
 \[
  \frac{1}{b\pi \left( \left(\frac{x-a}{b}\right)^2 + 1 \right)},\qquad x\in\RR;
 \]
 and the distribution of \ $\xi$ \ is denoted by \ $\PP_{C(a,b)}$.
\ In case of \ $a=0$ \ and \ $b=1$, \ $\PP_{C(a,b)}$ \ is called a standard Cauchy distribution.
Note that if \ $\xi$ \ has a standard Cauchy distribution, then \ $a\xi+b$ \ has distribution \ $\PP_{C(a,b)}$ \ for any \ $a\in\RR$ \ and $b>0$.

Extensions of Newton method for finding complex roots of polynomials has a long history, it goes back at least to the work of Cayley \cite{Cay1, Cay2}.
Let \ $a\in\RR$ \ and \ $b>0$.
The Newton iteration for the function \ $f:\RR\to\RR$, \ $f(x) := (x-a)^2 + b^2$, \ $x\in\RR$, \ takes the form
 \begin{align}\label{help1}
  \begin{split}
  x_{n+1} & = x_n - \frac{f(x_n)}{f'(x_n)}
           = x_n - \frac{(x_n-a)^2 + b^2}{2(x_n-a)}
           = \frac{b}{2}\left( \frac{x_n-a}{b} - \frac{1}{\frac{x_n-a}{b}} \right) + a
  \end{split}
 \end{align}
 for those \ $n\in\ZZ_+$ \ that satisfy \ $x_n\ne a$.
The two roots of the equation \ $f(z) = 0$, \ $z\in\CC$, \ are \ $a\pm \ii b$, \ and hence the
  perpendicular bisector of these two roots is the real axis.
By Cayley's theorem, the Newton iteration for \ $f$ \ started from a complex number having positive imaginary part
 converges to \ $a+\ii b$, \ while started from a complex number having negative imaginary part converges to \ $a-\ii b$.
Further, the Newton iteration for \ $f$ \ started from a real number do not converge (see, e.g., Devaney \cite[page 170]{Dev}),
 and below we recall an ergodic result (see Proposition \ref{Pro1}) which roughly speaking states that the arithmetic averages of appropriate functions
 of the Newton iterates in question converge.

For any \ $a\in\RR$ \ and \ $b>0$, \ let us introduce the function \ $\varphi_{a,b}:\RR\to\RR$,
 \begin{align}\label{Def_varphi}
   \varphi_{a,b}(x):=\begin{cases}
                       \frac{b}{2}\left( \frac{x-a}{b} - \frac{1}{\frac{x-a}{b}} \right) + a = \frac{x}{2} + \frac{a}{2} - \frac{\frac{b^2}{2}}{x-a}
                           & \text{if \ $x\ne a$,}\\
                       a   & \text{if \ $x =a$.}
                      \end{cases}
 \end{align}
Note that \ $\varphi_{a,b}$ \ is Borel measurable,  \ $\varphi_{a,b}(x) = b \varphi_{0,1}(\frac{x-a}{b}) + a$, \ $x\in\RR$, \ and
 one calls \ $\varphi_{a,b}$ \ an affine Boole transformation.
In the special case \ $(a,b)=(0,1)$, the map \ $\varphi_{0,1}$ \ is called Boole transformation.
The origin of the notion of Boole transformation is the following formula due to G. Boole \cite{Boo}:
 \[
   \int_\RR f(x)\,\dd x = \int_\RR f\left(x - \frac{1}{x}\right)\,\dd x
 \]
 for any integrable function \ $f:\RR\to\RR$.

It is known that \ $\varphi_{0,1}$ \ preserves the standard Cauchy distribution \ $\PP_{C(0,1)}$, \ see, e.g., Pitman and Williams \cite{PitWil}
 or Chin et al.\ \cite{ChiJunMar}.
It is also known that the Boole transformation \ $\varphi_{0,1}$ \ is ergodic with respect to the standard Cauchy distribution, see, e.g.,
  Prykarpatsky and Feldman \cite[Theorem 2.2]{PryFel}, Lee and Suriajaya \cite[Theorem 3.1]{LeeSur} or Chin et al.\ \cite[Lemma 9]{ChiJunMar}.

We also have that for all \ $a\in\RR$ \ and \ $b>0$, \ the map \ $\varphi_{a,b}$ \ preserves \ $\PP_{C(a,b)}$, \ i.e., the Cauchy distribution with
 parameter \ $(a,b)$, \ and is ergodic with respect to \ $\PP_{C(a,b)}$.
\ Indeed, since \ $\varphi_{0,1}$ \ preserves \ $\PP_{C(0,1)}$, \ for all \ $a\in\RR$ \ and $b>0$, \
 we have \ $\varphi_{a,b}$ \ preserves the Cauchy distribution with parameter \ $(a,b)$, \ i.e.,
 \ $\varphi_{a,b}(b\xi+a)\distre b\xi+a$, \ where \ $\xi$ \ is a random variable having distribution \ $\PP_{C(0,1)}$.
One can also easily check that \ $\varphi_{a,b}$ \ is ergodic with respect to the Cauchy distribution with parameter \ $(a,b)$.
Namely, we have to check that \ $\PP_{C(a,b)}(A) = 0$ \ or \ $\PP_{C(a,b)}(\RR\setminus A)=0$ \
 (equivalently \ $\PP_{C(a,b)}(A)\in\{0,1\}$) \ for any \ $A\in\cB(\RR)$ \ with \ $\varphi_{a,b}^{-1}(A) = A$.
\ Note that
 \begin{align*}
   \varphi_{a,b}^{-1}(A)
   & = \{x\in\RR : \varphi_{a,b}(x)\in A \}
     = \left\{x\in\RR : b\varphi_{0,1}\left(\frac{x-a}{b}\right)+a\in A \right\}\\
   & = \left\{x\in\RR : \varphi_{0,1}\left(\frac{x-a}{b}\right)\in \frac{A-a}{b} \right\}\\
   & = b\left\{\frac{x-a}{b} \in\RR : \varphi_{0,1}\left(\frac{x-a}{b}\right)\in \frac{A-a}{b} \right\} + a\\
   & = b\left\{ y \in\RR : \varphi_{0,1}(y)\in \frac{A-a}{b}\right \} + a
    = b \varphi_{0,1}^{-1}\left(\frac{A-a}{b}\right) + a.
 \end{align*}
Since \ $\varphi_{a,b}^{-1}(A)= A$, \ we have \ $\varphi^{-1}_{0,1}\left(\frac{A-a}{b}\right) = \frac{A-a}{b}$.
\ Hence, using that \ $\varphi_{0,1}$ \ is ergodic with respect to the standard Cauchy distribution  \ $\PP_{C(0,1)}$,
 \ we have \ $\PP_{C(0,1)}\left(\frac{A-a}{b}\right) \in\{0,1\}$, \ i.e., if \ $\xi$ \ has standard Cauchy distribution,
  \ then \ $\PP\left(\xi\in \frac{A-a}{b}\right)\in\{0,1\}$.
\ Consequently, using that \ $b\xi+a$ \ has Cauchy distribution with parameter \ $(a,b)$, \ we have
 \ $\PP_{C(a,b)}(A) = \PP(b\xi + a \in A) \in \{0,1\}$, \ as desired.

Using that \ $\varphi_{a,b}$ \ is measure preserving and ergodic with respect to the Cauchy distribution
 \ $\PP_{C(a,b)}$, \ by Birkhoff's ergodic theorem (see, e.g., Durrett \cite[Theorem 6.2.1]{Dur} or Einsiedler and Ward \cite[Theorem 2.30]{EinWar}),
 we have the following well-known result.

\begin{Pro}\label{Pro1}
For any \ $a\in\RR$ \ and \ $b>0$, \ let us consider the affine Boole transformation \ $\varphi_{a,b}$ \ given in \eqref{Def_varphi}.
Then for Lebesgue almost every \ $x\in\RR$ \ and for any Borel measurable function \ $f: \RR \to \RR$ \ such that
 \ $\int_\RR \vert f(u)\vert/( (u-a)^2/b^2+1)\,\dd u<\infty$, \ we have
 \[
    \frac{1}{n}\sum_{k=0}^{n-1} f(\varphi_{a,b}^{(k)}(x)) \to \int_\RR \frac{f(u)}{b\pi \left( \left(\frac{u-a}{b}\right)^2 + 1 \right)}\,\dd u
                                                               = \EE(f(\xi))
    \qquad \text{as \ $n\to\infty$,}
 \]
 where \ $\varphi_{a,b}^{(k)}$ \ denotes the \ $k$-fold iteration of \ $\varphi_{a,b}$ \ for \ $k\in\ZZ_+$ \
 with the convention \ $\varphi_{a,b}^{(0)}(x):=x$, \ $x\in\RR$, \ and \ $\xi$ \ is a random variable having Cauchy distribution with parameter \ $(a,b)$.
\end{Pro}

In the next remark, we recall some existing literature related to Proposition \ref{Pro1}.
We point out that Proposition \ref{Pro1} is well-known and authors reproved it several times independently of each other.

\begin{Rem}
(i). By Prykarpatsky and Feldman \cite[Lemma 2.1 and Theorem 2.2]{PryFel}, for any \ $\alpha\in \RR$ \ and \ $\beta>0$,
  \ the mapping \ $\varphi:\RR\to\RR$,
  \[
    \varphi(x):=\begin{cases}
             \frac{x}{2} + \alpha - \frac{\beta}{x-2\alpha}  & \text{if \ $x\ne 2\alpha$,}\\
              2\alpha  & \text{if \ $x=2\alpha$,}
          \end{cases}
  \]
 is measure preserving and ergodic with respect to the Cauchy distribution \ $\PP_{C(2\alpha,\sqrt{2\beta})}$.
One can easily check that \ $\varphi = \varphi_{2\alpha,\sqrt{2\beta}}$, \ so Proposition \ref{Pro1} can be applied to \ $\varphi$.

(ii). Proposition \ref{Pro1} is a special case of Theorem 4 in Ishitani and Ishitani \cite{IshIsh}.
In fact, Ishitani and Ishitani \cite[Example 1]{IshIsh} considered a generalized Boole transformation
 of the form $\RR\setminus\{0\}\ni x\mapsto \alpha x  -  \frac{\beta}{x}$, where $\alpha\in(0,1)$ and $\beta>0$, and
 derived an analogue of Proposition \ref{Pro1}.

(iii). Proposition \ref{Pro1} is a special case of the Theorem (Basic Ergodicity) in Umeno \cite[page 165]{Ume}.
In fact, Umeno \cite{Ume} proved similar results for a branch of transformations defined by the addition theorems for the tangent and cotangent functions
 such as for  \ $\RR\setminus\{\pm 1\}\ni x \mapsto 2x/(1-x^2)$ \ or \ $\RR\setminus\{\pm 1/\sqrt{3}\}\ni x \mapsto (3x-x^2)/(1-3x^2)$.

(iv). Proposition \ref{Pro1} is the same as Corollary 3.2 in Lee and Suriajaya \cite{LeeSur}.
Indeed, for \ $\alpha>0$ \ and \ $\beta\in\RR$, \ they consider the transformation \ $T_{\alpha,\beta}:\RR\to\RR$,
 \[
    T_{\alpha,\beta}(x):=\begin{cases}
                           \frac{\alpha}{2}\left(\frac{x+\beta}{\alpha} - \frac{\alpha}{x-\beta}\right)
                                   = \frac{x}{2} + \frac{\beta}{2} - \frac{\frac{\alpha^2}{2}}{x-\beta}
                                 & \text{if \ $x\ne \beta$,}\\
                           \beta & \text{if \ $x = \beta$,}
                          \end{cases}
 \]
 and hence \ $T_{\alpha,\beta} = \varphi_{\beta,\alpha}$.

(v).
By choosing $a=0$ and $b=1$ in Proposition \ref{Pro1}, we get back Proposition 6 in Chin et al.\ \cite{ChiJunMar}.
\proofend
\end{Rem}

In Appendix \ref{App_Examples} we give some applications of Proposition \ref{Pro1}.

Given a function \ $g:\RR\to\RR$, \ for each \ $k\in\ZZ_+$, \ let \ $g^{(k)}$ \ be the $k$-fold iteration
 of \ $g$ \ with the convention \ $g^{(0)}(x):=x$, \ $x\in\RR$.

\begin{Def}\label{Def_Lyapunov}
Let \ $g:\RR\to\RR$ \ be a function which is differentiable Lebesgue almost everywhere.
Let \ $x_0\in\RR$ \ be such that \ $g'(x_k)$ \ is well-defined and nonzero for each \ $k\in\ZZ_+$, \ where
 \ $x_k:=g^{(k)}(x_0)$, \ $k\in\ZZ_+$, \ and \ $g'(x_k)$ \ denotes the derivative of \ $g$ \ at \ $x_k$.
\ Then the Lyapunov exponent of the orbit \ $(g^{(k)}(x_0))_{k\in\ZZ_+}$ \ is defined by
 \[
   \lim_{n\to\infty} \frac{1}{n} \sum_{k=0}^{n-1} \ln(\vert g'(x_k)\vert),
 \]
 provided that the limit exists.
\end{Def}

Note that, under the conditions of Definition \ref{Def_Lyapunov}, if the limit \ $\lim_{n\to\infty} \frac{1}{n} \sum_{k=0}^{n-1} \ln(\vert g'(x_k)\vert)$
 \ exists, then, for each \ $m\in\NN$, \ we have the limit $\lim_{n\to\infty} \frac{1}{n} \sum_{k=0}^{n-1} \ln(\vert g'(x_{k+m})\vert)$ \ exists as well
 and the two limits are equal.
This can be interpreted in a way that the Lyapunov exponent is indeed a quantity corresponding
 to the orbit \ $(g^{(k)}(x_0))_{k\in\ZZ_+}$.

The next theorem states that the Lyapunov exponent of the affine Boole transformation \ $\varphi_{a,b}$ \
 (defined in \eqref{Def_varphi}) is \ $\ln(2)$, \ which follows, e.g., from Umeno and Okubo \cite[Theorem 2]{UmeOku}.

\begin{Thm}\label{Thm1}
For all \ $a\in\RR$ \ and \ $b>0$, \ the Lyapunov exponent of the affine Boole transformation
 \ $\varphi_{a,b}$ \ defined in \eqref{Def_varphi} is \ $\ln(2)$, \ more precisely, we have
 \[
    \lim_{n\to\infty}\frac{1}{n}\sum_{k=0}^{n-1} \ln(\vert \varphi_{a,b}'(x_k)\vert) = \ln(2)
 \]
 for Lebesgue almost every \ $x_0\in\RR$, \ where \ $x_k:=\varphi_{a,b}^{(k)}(x_0)$, \ $k\in\ZZ_+$,
 \ and \ $\varphi_{a,b}^{(k)}$ \ denotes the \ $k$-fold iteration of \ $\varphi_{a,b}$ \ for \ $k\in\ZZ_+$ \
 with the convention \ $\varphi_{a,b}^{(0)}(x) := x$, \ $x\in\RR$.
\end{Thm}

In the present paper we give an elementary and alternative proof of Theorem \ref{Thm1}
 compared to the existing ones available in the literature.
Our proof is based on Proposition \ref{Pro1} and we calculate the corresponding integral appearing in Proposition \ref{Pro1}
 using the method of differentiation under the integral sign.
It turns out that we need to calculate the integral \ $\int_0^\infty \frac{\ln(1+x^{-2})}{\pi(1+x^2)}\,\dd x$, \
 which can be transformed into several other forms such as \ $-\frac{2}{\pi} \int_0^{\pi/2} \ln(\sin(x))\,\dd x$, \
 $\frac{2}{\pi} \int_0^{\pi/2} \frac{x}{\tan(x)}\,\dd x$ \ or \ $\frac{2}{\pi}\int_0^\infty \frac{\arctan(x)}{x(1+x^2)}\,\dd x$, \
 see, e.g., Blackmore et al.\ \cite[formulae (44) and (65)]{BlaBalKycPry}.
As we already mentioned, it is well-known that the common value of these integrals is $\ln(2)$,
 which now we calculate using the method of differentiation under the integral sign.

For historical fidelity, we note that Umeno and Okubo \cite[Theorem 2]{UmeOku} calculated the Lyapunov exponent of a generalized Boole
 transformation $\RR\setminus\{0\}\ni x \mapsto \alpha x - \frac{\beta}{x}$, where $\alpha\in(0,1)$ and $\beta>0$ using Proposition \ref{Pro1}.
However, their proof technique is completely different, they used complex analysis for the calculation of the corresponding integral
 appearing in Proposition \ref{Pro1}.
We also mention that, for all $a\in\RR$ and $b>0$, the metric (also called measure-theoretic) entropy on the probability space
 \ $(\RR,\cB(\RR),\PP_{C(a,b)})$ \ of the affine Boole transformation \ $\varphi_{a,b}$ \ defined in \eqref{Def_varphi} is \ $\ln(2)$
 \ as well, see, e.g., Blackmore et al.\ \cite[page 13]{BlaBalKycPry}.

The remaining part of the paper is structured as follows.
Section \ref{section_proof} contains a proof of Theorem \ref{Thm1} using the method of differentiation under the integral sign.
In Appendix \ref{App_Klenke}, we recall results on the continuity and differentiability of parametric integrals that we use in the proof of Theorem \ref{Thm1}.
In Appendix \ref{App_Examples}, we give some further applications of Proposition \ref{Pro1} for some particular functions \ $f$.

\section{Proof of Theorem \ref{Thm1} using the method of differentiation under the integral sign}
\label{section_proof}

First, we prove Theorem \ref{Thm1} in case of \ $(a,b)=(0,1)$.
\ Note that the function \ $\varphi_{0,1}$ \ is differentiable on \ $\RR\setminus\{0\}$ \
 with $\varphi_{0,1}'(x) = \frac{1}{2}\big(1+\frac{1}{x^2}\big)\ne 0$, $x\in\RR\setminus\{0\}$,
 \ but it is not differentiable at \ $0$ \ following from
 \[
   \lim_{x\to 0} \frac{\varphi_{0,1}(x) - \varphi_{0,1}(0)}{x-0}
                 = \lim_{x\to 0} \frac{ \frac{1}{2} \left(x-\frac{1}{x}\right) - 0}{x-0}
                 = \lim_{x\to 0} \frac{1}{2}\left(1-\frac{1}{x^2}\right)
                 =-\infty.
 \]
However, we can check that the sequence \ $(\varphi_{0,1}'(x_k))_{k\in\ZZ_+} = (\varphi_{0,1}'(\varphi_{0,1}^{(k)}(x_0)))_{k\in\ZZ_+}$ \ is well-defined
 and nonzero for all \ $x_0\in\RR\setminus N$, \ where \ $N$ \ is a countable subset of \ $\RR$.
\ For this it is enough to check that the set
 \begin{align}\label{help_countable}
   \Big\{x_0\in\RR : \varphi_{0,1}^{(k)}(x_0) = 0 \;\; \text{for some \ $k\in\ZZ_+$}\Big\}
 \end{align}
 is countable.
This set can be written in the form \ $\bigcup_{k=1}^\infty A_k$, \ where
 \[
  A_k:=\big\{x_0\in\RR : \varphi_{0,1}^{(\ell)}(x_0) = 0 \;\; \text{for some \ $\ell\in\{0,1,\ldots,k\}$}\big\}, \qquad k\in\NN.
 \]
Using that the union of countably many finite sets is countable, it is enough to verify that \ $A_k$ \ is finite for each \ $k\in\NN$.
\ We prove it by induction.
We have \ $\varphi_{0,1}^{(0)}(x)=0$ \ if and only if \ $x=0$.
\ Further, \ $\varphi_{0,1}(x)=0$ \ if and only if \ $x=0$ \ or \ $\frac{1}{2}\left(x-\frac{1}{x}\right)=0$, \ $x\ne 0$, \ which holds
 if and only if \ $x=\pm 1$, \  and hence \ $A_1=\{-1, 0, 1\}$ \ being a finite set.
Suppose that the sets \ $A_1,\ldots,A_k$ \ are finite.
Then \ $A_{k+1} = A_k\cup\{x_0\in\RR : \varphi_{0,1}^{(k+1)}(x_0)=0 \}$, \ and the function
 \[
 \RR\setminus \bigcup_{\ell=1}^k A_\ell \ni x \mapsto \varphi_{0,1}^{(k+1)}(x)= \frac{1}{2}\left( \varphi_{0,1}^{(k)}(x)  - \frac{1}{\varphi_{0,1}^{(k)}(x)}\right)
 \]
 is a well-defined rational function such that its numerator is a polynomial of degree \ $2^{k+1}$ \ and its denominator is a polynomial
 of degree \ $2^{k+1}-1$.
\ Consequently, the equation \ $\varphi_{0,1}^{(k+1)}(x) = 0$, \ $x\in\RR\setminus \bigcup_{\ell=1}^k A_\ell$, \ has at most
 \ $2^{k+1}$ \ (i.e., finitely many) solutions due to the fundamental theorem of algebra.
Hence we have \ $A_{k+1}$ \ is finite, as desired.

Let us apply Proposition \ref{Pro1} with \ $a=0$, \ $b=1$, \ and the Borel measurable function \ $f:\RR\to\RR$,
 \begin{align}\label{help_f}
    f(x):=
           \begin{cases}
               \ln(\vert \varphi_{0,1}'(x)\vert) = \ln\left(\frac{1}{2}\left(1+\frac{1}{x^2}\right)\right) & \text{if \ $x\in\RR\setminus\{0\}$},\\
               0 & \text{if \ $x=0$.}
           \end{cases}
 \end{align}
Note that we cannot define \ $f$ \ to be \ $\ln(\vert \varphi_{0,1}'(x)\vert)$ \ for all \ $x\in\RR$, \
 since, as we noted, the function \ $\varphi_{0,1}$ \ is not differentiable at \ $0$.
\ Using that
 \begin{align*}
    \int_\RR \frac{\vert f(x)\vert }{\pi(1+x^2)}\,\dd x
     & = \int_{(-\infty,0)} \frac{\vert \ln\left(\frac{1}{2}\left(1+\frac{1}{x^2}\right)\right)\vert }{\pi(1+x^2)}\,\dd x
         + \int_{(0,\infty)} \frac{\vert \ln\left(\frac{1}{2}\left(1+\frac{1}{x^2}\right)\right)\vert}{\pi(1+x^2)}\,\dd x \\
     & = 2\int_{(0,\infty)} \frac{\vert \ln\left(\frac{1}{2}\left(1+\frac{1}{x^2}\right)\right) \vert }{\pi(1+x^2)}\,\dd x \\
     &\leq 2\ln(2) \int_{(0,\infty)} \frac{1}{\pi(1+x^2)}\,\dd x
           + 2\int_{(0,\infty)} \frac{\ln\left(1+\frac{1}{x^2}\right)}{\pi(1+x^2)}\,\dd x \\
     & = \ln(2) + 2\int_{(0,\infty)} \frac{\ln\left(1+\frac{1}{x^2}\right)}{\pi(1+x^2)}\,\dd x ,
 \end{align*}
 in order to have right to apply Proposition \ref{Pro1} we need to check that
 \begin{align}\label{help2_Feynman}
 \int_{(0,\infty)} \frac{\ln\left(1+\frac{1}{x^2}\right)}{\pi(1+x^2)}\,\dd x < \infty.
 \end{align}

We calculate the integral in \eqref{help2_Feynman} using the method of differentiation under the integral sign.
Let \ $G:[0,1]\to\RR$ \ be given by
 \[
   G(t):=  \int_{(0,\infty)} \frac{\ln\left(1+\frac{t}{x^2}\right)}{\pi(1+x^2)}\,\dd x, \qquad t\in[0,1].
 \]
Then we need to calculate \ $G(1)$, \ and note that \ $G(0)=0$.
\ Let \ $\vare\in(0,1)$ \ be arbitrary, and let us apply Theorem \ref{Thm_differentiation} with the following choices
 \begin{align}\label{help_ThmA2_choices}
  \begin{split}
    &I:=(\vare,1),\qquad X:=(0,\infty), \qquad \cA:=\cB((0,\infty)),\qquad \mu:=\text{Lebesgue measure on \ $\cA$},\\
    &g(x,t):=\frac{\ln\left(1+\frac{t}{x^2}\right)}{\pi(1+x^2)},\quad x\in(0,\infty),\;\; t\in(\vare,1).
   \end{split}
 \end{align}
Next we check the conditions of Theorem \ref{Thm_differentiation}.
For any \ $t\in(\vare,1)$, \ by partial integration,
 \begin{align}\label{help4}
 \begin{split}
   0 & \leq \int_{(0,\infty)} \frac{\ln\left(1+\frac{t}{x^2}\right)}{\pi(1+x^2)} \,\dd x
      \leq \frac{1}{\pi} \int_{(0,\infty)} \ln\left(1+\frac{t}{x^2}\right) \,\dd x \\
   & = \lim_{x\to\infty} \frac{x}{\pi}\ln\left(1+\frac{t}{x^2}\right)
          -  \lim_{x\downarrow 0} \frac{x}{\pi}\ln\left(1+\frac{t}{x^2}\right)
          + 2t\int_{(0,\infty)} \frac{1}{\pi(x^2 + t)}\,\dd x\\
   & =  2\int_{(0,\infty)} \frac{1}{\pi((x/\sqrt{t})^2 + 1)}\,\dd x
      = 2 \sqrt{t} \int_{(0,\infty)} \frac{1}{\pi(y^2 + 1)}\,\dd y
     =   \sqrt{t} <\infty,
  \end{split}
 \end{align}
 since, by L'Hospital's rule,
 \[
  \lim_{x\to\infty} x\ln\left(1+\frac{t}{x^2}\right)
     =  \lim_{x\to\infty} \frac{\ln\left(1+\frac{t}{x^2}\right)}{\frac{1}{x}}
     = \lim_{x\to\infty}  \frac{2t}{x+\frac{t}{x}} = 0
 \]
 and
 \[
  \lim_{x\downarrow 0} x\ln\left(1+\frac{t}{x^2}\right)
     = \lim_{x\downarrow 0}  \frac{2t}{x+\frac{t}{x}} = 0.
 \]
Hence condition (i) of Theorem \ref{Thm_differentiation} is satisfied.

For all \ $x\in(0,\infty)$, \ the map \ $(\vare,1)\ni t \mapsto \frac{\ln\left(1+\frac{t}{x^2}\right)}{\pi(1+x^2)}$ \ is differentiable
 with derivative
 \[
     (\vare,1)\ni t\mapsto \frac{1}{\pi (1+x^2)(t+x^2)}=\partial_2g(x,t),
 \]
 and hence condition (ii) of Theorem \ref{Thm_differentiation} is satisfied.

For all \ $t\in(\vare,1)$, \ we have
 \[
    \vert \partial_2g(x,t)\vert \leq \frac{1}{\pi(1+x^2)t}\leq \frac{1}{\pi\vare(1+x^2)},\qquad x\in(0,\infty),
 \]
 where \ $\frac{1}{\pi\vare(1+x^2)}$, \ $x\in(0,\infty)$, \ is integrable on \ $(0,\infty)$ \ with respect to the Lebesgue measure,
 so condition (iii) of Theorem \ref{Thm_differentiation} is satisfied.

Thus we can apply Theorem \ref{Thm_differentiation} and we have
 \[
   G'(t) = \int_{(0,\infty)} \frac{1}{\pi(1+x^2)(t+ x^2)} \,\dd x, \qquad t\in(\vare,1).
 \]
Using that
 \[
   \frac{1}{(1+x^2)(t+ x^2)}  = \frac{1}{1-t}\left( - \frac{1}{1+x^2} + \frac{1}{t+x^2}\right),
   \qquad t\in(\vare,1), \; x\in\RR,
 \]
 we have
 \begin{align*}
  \int_{(0,\infty)} \frac{1}{\pi(1+x^2)(t+ x^2)}  \,\dd x
   & = \frac{1}{\pi(1-t)} \int_{(0,\infty)} \left(-\frac{1}{1+x^2} + \frac{1}{t+x^2} \right)\,\dd x \\
   & = \frac{1}{\pi(1-t)} \left(-\frac{\pi}{2}\right)
      + \frac{1}{\pi t(1-t)} \int_{(0,\infty)}\frac{1}{1+\Big(\frac{x}{\sqrt{t}}\Big)^2}\,\dd x
 \end{align*}
 \begin{align*}
   & = -\frac{1}{2(1-t)}
       + \frac{1}{\pi t(1-t)} \int_{(0,\infty)} \frac{1}{1+y^2}\sqrt{t}\,\dd y\\
   & = - \frac{1}{2(1-t)} +  \frac{1}{2\sqrt{t}(1-t)}
    = \frac{1}{2\sqrt{t}(1+\sqrt{t})}, \qquad t\in(\vare,1).
 \end{align*}
Consequently,
 \[
   G'(t) = \frac{1}{2\sqrt{t}(1+\sqrt{t})}, \qquad t\in(\vare,1),
 \]
 and hence
 \[
   G(1) - G(\vare) = \int_\vare^1 \frac{1}{2\sqrt{t}(1+\sqrt{t})} \,\dd t.
 \]
Then, by substituting \ $\sqrt{t}=z$, we have
 \begin{align}\label{help5}
  \begin{split}
   G(1) - G(\vare)
        =\int_{\sqrt{\vare}}^1 \frac{1}{2z(1+z)}2z\,\dd z
        = \int_{\sqrt{\vare}}^1 \frac{1}{1+z}\,\dd z
        = \ln(2) - \ln(1+\sqrt{\vare}).
  \end{split}
 \end{align}

Using Theorem \ref{Thm_continuity} we check that the function \ $G$ \ is continuous at \ $0$.
\ Let us apply Theorem \ref{Thm_continuity} with the following choices: \ $E:=[0,1]$, \ $d$ \ is the usual
 Euclidean metric, \ $t_0:=0$, \ $X$, \ $\cA$, \ and \ $\mu$ \ are given in \eqref{help_ThmA2_choices},
 and \ $g$ \ is also given in \eqref{help_ThmA2_choices} with the extension of its domain to \ $(0,\infty)\times[0,1]$.
\ Then condition (i) of Theorem \ref{Thm_continuity} holds (it follows from \eqref{help4} in case of \ $t\in(0,1]$, \
  and from \ $g(x,0) = 0$, \ $x\in(0,\infty)$, \ in case of \ $t=0$).
\ Condition (ii) of Theorem \ref{Thm_continuity} readily holds.
Condition (iii) of Theorem \ref{Thm_continuity} holds as well, since
 \begin{align*}
   \left\vert  \frac{\ln\left(1+\frac{t}{x^2}\right)}{\pi(1+x^2)}  \right\vert
    \leq \ln\left(1+\frac{t}{x^2}\right) \leq \ln\left(1+\frac{1}{x^2}\right),
    \qquad t\in[0,1],\;\; x\in(0,\infty),
 \end{align*}
 and, similarly to \eqref{help4}, we have
 \[
     \int_{(0,\infty)} \ln\left(1+\frac{1}{x^2}\right) \,\dd x =\pi <\infty.
 \]
Thus we can apply Theorem \ref{Thm_continuity} and we have the continuity of \ $G$ \ at \ $0$.

Consequently, by taking the limit as \ $\vare\downarrow 0$ \ in \eqref{help5}, and using \ $G(0)=0$, \ we have
 \begin{align}\label{help6}
    G(1) = G(1) - G(0) = \ln(2) - \ln(1) = \ln(2),
 \end{align}
 and especailly, we have \eqref{help2_Feynman}.

Hence we can apply Proposition \ref{Pro1} with \ $a=0$, \ $b=1$ \ and the function \ $f$ \ given in \eqref{help_f}, and we get
 \[
    \frac{1}{n}\sum_{k=0}^{n-1} f(\varphi_{0,1}^{(k)}(x_0)) \to \int_\RR \frac{f(x)}{\pi(1+x^2)}\,\dd x
    \qquad \text{as \ $n\to\infty$}
 \]
 for Lebesgue almost every \ $x_0\in\RR$.
\ Since, by \eqref{help6},
 \begin{align}\label{help7}
  \begin{split}
   \int_\RR \frac{f(x)}{\pi(1+x^2)}\,\dd x
     &= 2\int_{(0,\infty)} \frac{\ln\left(\frac{1}{2}\left(1+\frac{1}{x^2}\right)\right)}{\pi(1+x^2)}\,\dd x
      = 2\left( -\ln(2) \int_0^\infty \frac{1}{\pi(1+x^2)}\,\dd x + G(1)\right) \\
     &= 2\left(\frac{-\ln(2)}{2} + G(1)\right) = \ln(2),
  \end{split}
 \end{align}
 we have
 \[
    \frac{1}{n}\sum_{k=0}^{n-1} f(\varphi_{0,1}^{(k)}(x_0)) \to \ln(2) \qquad \text{as \ $n\to\infty$}
 \]
 for Lebesgue almost every \ $x_0\in\RR$.
\ Using that the set in \eqref{help_countable} is countable, we have
 \begin{align}\label{help3}
    \frac{1}{n}\sum_{k=0}^{n-1} \ln(\vert \varphi'_{0,1}(x_k)\vert) \to \ln(2) \qquad \text{as \ $n\to\infty$}
 \end{align}
 for Lebesgue almost every \ $x_0\in\RR$, \ where \ $x_k:=\varphi_{0,1}^{(k)}(x_0)$, \ $k\in\ZZ_+$, \  as desired.

Now we turn to prove Theorem \ref{Thm1} in the general case \ $(a,b)\in\RR\times(0,\infty)$.
By induction, one can check that
 \[
  \varphi_{a,b}^{(k)}(y_0) = b \varphi_{0,1}^{(k)}\left(\frac{y_0-a}{b}\right)+a,\qquad k\in\ZZ_+,
 \]
 and further, we have that \ $\varphi_{a,b}'(x) = \varphi_{0,1}'\big(\frac{x-a}{b}\big)$ \ for \ $x\in\RR\setminus\{a\}$.
\ Hence, using that the set in \eqref{help_countable} is countable, we have that the set
 \begin{align*}
   \Big\{y_0\in\RR : \varphi_{a,b}^{(k)}(y_0) = a \;\; \text{for some \ $k\in\ZZ_+$}\Big\}
    & = \left\{y_0\in\RR : \varphi_{0,1}^{(k)}\left(\frac{y_0-a}{b}\right) = 0 \;\; \text{for some \ $k\in\ZZ_+$}\right\}\\
    & = b\Big\{x_0\in\RR : \varphi_{0,1}^{(k)}(x_0) = 0 \;\; \text{for some \ $k\in\ZZ_+$}\Big\} + a
 \end{align*}
 is countable as well.
Therefore, by \eqref{help3}, we get
 \begin{align*}
   \frac{1}{n}\sum_{k=0}^{n-1} \ln( \vert \varphi_{a,b}'(y_k)\vert )
    & = \frac{1}{n}\sum_{k=0}^{n-1} \ln\left( \left\vert \varphi_{0,1}'\left(\frac{y_k-a}{b}\right)\right\vert \right)
     = \frac{1}{n}\sum_{k=0}^{n-1} \ln\left( \left\vert \varphi_{0,1}'\left(  \varphi_{0,1}^{(k)}\left(\frac{y_0-a}{b}\right) \right)\right\vert \right)\\
   & \to \ln(2)
 \end{align*}
 as \ $n\to\infty$ \ for Lebesgue almost every \ $y_0\in\RR$, \ where \ $y_k:=\varphi_{a,b}^{(k)}(y_0)$, \ $k\in\ZZ_+$, \ as desired.
\proofend

\appendix

\vspace*{5mm}

\noindent{\bf\Large Appendix}

\section{Continuity and differentiation of parametric integrals}\label{App_Klenke}

We recall two results on the continuity and differentiability of parametric integrals of functions of two variables where
 the integration is taken with respect to one of the variables, see, e.g., Klenke \cite[Theorems 6.27 and 6.28]{Kle}.
Given a measure space \ $(X,\cA,\mu)$, \ a function \ $f:X\to\RR$ \ is said to be in \ $L^1(X,\cA,\mu)$ \ (for short in \ $L^1(\mu)$)
 \ if it is \ $(\cA,\cB(\RR))$-measurable and \ $\int_X \vert f(x)\vert\,\mu(\dd x)<\infty$.

\begin{Thm}[continuity of parametric integrals]\label{Thm_continuity}
Let \ $(X,\cA,\mu)$ \ be a \ $\sigma$-finite measure space, \ $(E,d)$ \ be a metric space, \ $t_0\in E$, \ and
 let \ $g:X\times E\to\RR$ \ be a map with the following properties:
 \begin{itemize}
   \item[(i)] for any \ $t\in E$, \ the map \ $X\ni x\mapsto g(x,t)$ \ is in \ $L^1(\mu)$,
   \item[(ii)] for \ $\mu$-almost every \ $x\in X$, \ the map \ $E\ni t\mapsto g(x,t)$ \ is continuous at \ $t_0$,
   \item[(iii)] there is a map \ $h:X\to\RR$ \ such that \ $h\geq 0$, \ $h\in L^1(\mu)$ \ and
                \[
                 \vert g(x,t)\vert \leq h(x)\qquad \text{$\mu$-a.e. \ $x\in X$ \ for all \ $t\in E$.}
                \]
 \end{itemize}
Then the map \ $E\ni t\mapsto \int_X g(x,t)\,\mu(\dd x)$ \ is continuous at \ $t_0$.
\end{Thm}

\begin{Thm}[differentiability of parametric integrals]\label{Thm_differentiation}
Let \ $(X,\cA,\mu)$ \ be a \ $\sigma$-finite measure space, \ $I$ \ be a nontrivial open interval (having at least two different points),
 and let \ $g:X\times I\to\RR$ \ be a map with the following properties:
 \begin{itemize}
   \item[(i)] for any \ $t\in I$, \ the map \ $X\ni x\mapsto g(x,t)$ \ is in \ $L^1(\mu)$,
   \item[(ii)] for \ $\mu$-almost every \ $x\in X$, \ the map \ $I\ni t\mapsto g(x,t)$ \ is differentiable with derivative denoted by \ $I\ni t\mapsto \partial_2 g(x,t)$,
   \item[(iii)] there is a map \ $h:X\to\RR$ \ such that \ $h\geq 0$, \ $h\in L^1(\mu)$ \ and
                \[
                 \vert \partial_2 g(x,t)\vert \leq h(x)\qquad \text{$\mu$-a.e. \ $x\in X$ \ for all \ $t\in I$.}
                \]
 \end{itemize}
Then the map
 \ $G: I\to\RR$, \ $G(t):=\int_X g(x,t)\,\mu(\dd x)$, \ $t\in I$, \ is differentiable with derivative
 \[
   G'(t) = \int_X \partial_2 g(x,t)\,\mu(\dd x),\qquad t\in I.
 \]
\end{Thm}

\section{Some applications of Proposition \ref{Pro1}}\label{App_Examples}

We give some applications of Proposition \ref{Pro1}.
Recall that for all \ $a\in\RR$ \ and \ $b>0$, \ $\varphi_{a,b}^{(k)}$ \ denotes the \ $k$-fold iteration
 of \ $\varphi_{a,b}$ \ for \ $k\in\ZZ_+$ \ with the convention \ $\varphi_{a,b}^{(0)}(x) := x$, \ $x\in\RR$.
\begin{Ex}
Let \ $a\in\RR$ \ and \ $b>0$.
     \ If \ $f(u):= b\pi \left( \left(\frac{u-a}{b}\right)^2 + 1 \right)\frac{1}{\sqrt{2\pi}} \ee^{-\frac{u^2}{2}}$, \ $u\in\RR$, \ then
    \[
      \int_\RR \frac{\vert f(u)\vert}{b\pi \left( \left(\frac{u-a}{b}\right)^2 + 1 \right)}\,\dd u
       = \int_\RR \frac{f(u)}{b\pi \left( \left(\frac{u-a}{b}\right)^2 + 1 \right)}\,\dd u
       = \int_\RR \frac{1}{\sqrt{2\pi}} \ee^{-\frac{u^2}{2}}\,\dd u
       = 1,
    \]
    so, by Proposition \ref{Pro1}, for Lebesgue almost every \ $x\in\RR$, \ we have
    \[
       \frac{1}{n}\sum_{k=0}^{n-1} f(\varphi_{a,b}^{(k)}(x)) \to 1 \qquad \text{as \ $n\to\infty$.}
    \vspace*{-7mm}
    \]
\proofend
\end{Ex}

\begin{Ex}
Let \ $a\in\RR$ \ and \ $b:=1$.
   If \ $f(u):=\pi ((u-a)^2 + 1) u \frac{1}{\sqrt{2\pi}} \ee^{-\frac{(u-a)^2}{2}}$, \ $u\in\RR$, \ then
   using \ $\vert u\vert\leq u^2+1$, \ $u\in\RR$, \ we have
   \begin{align*}
     \int_\RR \frac{\vert f(u)\vert}{b\pi \left( \left(\frac{u-a}{b}\right)^2 + 1 \right)}\,\dd u
      & = \int_\RR \frac{\vert f(u)\vert }{\pi((u-a)^2 +1)}\,\dd u
      = \int_\RR \vert u\vert \frac{1}{\sqrt{2\pi}} \ee^{-\frac{(u-a)^2}{2}} \,\dd u \\
     & \leq \EE( (\zeta + a)^2 ) +1 = a^2 + 2,
   \end{align*}
   where \ $\zeta$ \ is a standard normally distributed random variable, and
   \[
     \int_\RR \frac{f(u)}{b\pi \left( \left(\frac{u-a}{b}\right)^2 + 1 \right)}\,\dd u
      = \int_\RR \frac{f(u)}{\pi((u-a)^2 +1)}\,\dd u = \int_\RR u \frac{1}{\sqrt{2\pi}} \ee^{-\frac{(u-a)^2}{2}} \,\dd u = a.
   \]
   Hence, by Proposition \ref{Pro1}, for Lebesgue almost every \ $x\in\RR$, \ we have
    \[
       \frac{1}{n}\sum_{k=0}^{n-1} f(\varphi_{a,1}^{(k)}(x)) \to a \qquad \text{as \ $n\to\infty$.}
    \vspace*{-7mm}
    \]
\proofend
\end{Ex}

\begin{Ex}
Let \ $a\in\RR$, \ $b:=1$, \ and \ $\eta$ \ be an absolutely continuous random variable with a density function \ $h_\eta$ \ such that
 \ $\EE(\eta^2)<\infty$. \
If \ $f(u):=\pi(u^2+1)h_\eta(u)$, \ $u\in\RR$, \
         then
         \begin{align*}
           \int_\RR \frac{\vert f(u)\vert}{b\pi\left(\left(\frac{u-a}{b}\right)^2 +1\right)}\,\dd u
            & = \int_\RR \frac{f(u)}{\pi((u-a)^2 +1)}\,\dd u = \int_\RR\frac{1+u^2}{1+(u-a)^2} h_\eta(u) \,\dd u\\
            & = \EE\left(\frac{1+\eta^2}{1+(\eta-a)^2}\right)
            \leq \EE(1+\eta^2) <\infty.
         \end{align*}
         Hence, by Proposition \ref{Pro1}, for Lebesgue almost every \ $x\in\RR$, \ we have
        \[
           \frac{1}{n}\sum_{k=0}^{n-1} f(\varphi_{a,1}^{(k)}(x)) \to \EE\left(\frac{1+\eta^2}{1+(\eta-a)^2}\right) \qquad \text{as \ $n\to\infty$.}
        \]
        Here the function \ $\RR\ni a \mapsto \frac{1}{\pi\EE(1+\eta^2)}\EE\left(\frac{1+\eta^2}{1+(\eta-a)^2}\right)$ \ is a bounded density function,
        since it is Borel measurable, non-negative and, by Fubini's theorem,
        \begin{align*}
          &\int_{-\infty}^\infty \EE\left(\frac{1+\eta^2}{1+(\eta-a)^2}\right) \,\dd a
             = \int_{-\infty}^\infty \left(\int_\Omega  \frac{1+\eta(\omega)^2}{1+(\eta(\omega)-a)^2} \,\PP(\dd\omega) \right) \,\dd a\\
          & = \int_\Omega \left(\int_{-\infty}^\infty   \frac{1+\eta(\omega)^2}{1+(\eta(\omega)-a)^2} \, \,\dd a \right) \PP(\dd\omega) \\
          &= \int_\Omega (1+\eta(\omega)^2)\Big(\lim_{a\to\infty} \Big(-\arctan(\eta(\omega)-a)\Big) + \lim_{a\to-\infty} \Big(\arctan(\eta(\omega)-a)\Big) \Big) \, \PP(\dd\omega) \\
          &=\pi \int_\Omega (1+\eta(\omega)^2) \, \PP(\dd\omega) =\pi\EE(1+\eta^2).
        \end{align*}
        Further,
        \[
          \frac{1}{\pi\EE(1+\eta^2)} \EE\left(\frac{1+\eta^2}{1+(\eta-a)^2}\right) \leq \frac{1}{\pi},\qquad a\in\RR.
        \vspace*{-3mm}
        \]
\proofend
\end{Ex}

\section*{Acknowledgements}
I would like to thank Csaba Nosz\'aly for calling my attention to a possible application of Proposition \ref{Pro1} for
 calculation of Lyapunov exponents.
I would like to thank the referee for the comments that helped me improve the paper.

\section*{Declaration of competing interest}
The author declares that he has no known competing financial interests or personal relationships
that could have appeared to influence the work reported in this paper.

\end{document}